\documentclass[10pt,onecolumn]{IEEEtran}
\usepackage{url}
\newcommand{\citep}[1]{\cite{#1}}
\usepackage{setspace}

\usepackage{verbatim}
\usepackage{algorithm}
\usepackage{bm}
\usepackage{algpseudocode}
\usepackage[all]{xy}
\usepackage{url}
\usepackage{amsmath,amsthm,amssymb}
\usepackage{graphicx}
\usepackage{subcaption}

\newtheorem{thm}{Theorem}[section]
\newtheorem{theorem}{Theorem}[section]

\newtheorem{corollary}{Corollary}

\newtheorem{lemma}{lemma}
\theoremstyle{definition}
\newtheorem{definition}{Definition}

\theoremstyle{remark}
\newtheorem{remark}{Remark} 

\newcommand{\D}{\mathcal{D}}

\newcommand{\brs}[1]{\left[{#1}\right]}

\newcommand{\opt}[1]{{#1}^\ast}

\newcommand{\expb}[1]{\exp\left({#1}\right)}
\newcommand{\logb}[1]{\log\left({#1}\right)}

\newcommand{\powb}[2]{{\left({#1}\right)}^{#2}}

\newcommand{\tran}[1]{{#1}^T}

\newcommand{\ic}{\mathbf{j}}

\newcommand{\detb}[1]{\det\left({#1}\right)}
\newcommand{\tranb}[1]{{\left({#1}\right)}^T}
\newcommand{\br}[1]{\left({#1}\right)}

\newcommand{\inv}[1]{{\left(#1\right)}^{-1}}
\newcommand{\inve}[1]{{#1}^{-1}}

\newcommand{\tr}[1]{\mathrm{tr}\left( {#1} \right)}

\newcommand{\costt}{\ell}
\newcommand{\cost}{\costt_t}

\newcommand{\R}{\mathbb{R}}

\newcommand{\C}{\mathcal{C}}

\newcommand{\FHinf}[1]{q_\infty}
\newcommand{\FHtwo}[1]{q_2}
\newcommand{\FHone}[1]{q_1}

\newcommand{\E}{\mathcal{E}}

\newcommand{\G}{\mathrm{pv}}

\newcommand{\Vs}{\mathcal{V}}

\newcommand{\Lo}{\mathrm{pq}}

\newcommand{\Vset}{v}

\newcommand{\PQ}{$(P,Q)$~}
\newcommand{\PV}{$(P,V)$~}

\newcommand{\Int}{\mathrm{Int}}

\newcommand{\Sb}{\mathcal{S}}

\graphicspath{ {Figures/} }

\title{Convexity of Energy-Like Functions: Theoretical Results and Applications to Power System Operations}
\author{Krishnamurthy~Dvijotham,~\IEEEmembership{Member,~IEEE,}
        Steven~Low,~\IEEEmembership{Fellow,IEEE,}
        and~Michael~Chertkov,~\IEEEmembership{Member,~IEEE}
\thanks{K. Dvijotham and S. Low are with the Department
Computing and Mathematical Sciences, California Institute of Technology, CA, 91125 USA e-mail: dvij@caltech.edu,slow@caltech.edu.}
\thanks{M. Chertkov is with the T-Division, Los Alamos National Laboratory, Los Alamos, NM, USA 87544.}}
\begin{document}

\maketitle

\begin{abstract}
Power systems are undergoing unprecedented transformations with the incorporation of larger amounts of renewable energy sources, distributed generation and demand response. All these changes, while potentially making power grids more responsive, efficient and resilient, also pose significant implementation challenges. In particular, operating the new power grid will require new tools and algorithms capable of predicting if the current state of the system is operationally safe. In this paper we study and generalize the so-called energy function as a tool to design algorithms to test if a high-voltage power transmission system is within the allowed operational limits. In the past the energy function technique was utilized primarily to access the power system transient stability. In this manuscript, we take a new look at energy functions and focus on an aspect that has previously received little attention: \emph{Convexity}. We characterize the domain of voltage magnitudes and phases within which the energy function is convex. We show that the domain of the energy function convexity is sufficiently large to include most operationally relevant and practically interesting cases. We show how the energy function convexity can be used to analyze power flow equations, e.g. to certify solution uniqueness or non-existence  within the domain of convexity.  This and other useful features of the generalized energy function are described and illustrated on IEEE 14 and 118 bus models.
\end{abstract} 
\section{Introduction}

Power systems are experiencing revolutionary changes. Integration of renewable generation, distributed generation, smart metering, direct or price-based load-control capabilities have contributed to the revolution, but are also posing significant operational challenges by making the power system inherently stochastic and inhomogeneous. With these changes, the system operators will no longer have the luxury of large positive and negative reserves. Moreover, operating the future power grid will require developing new computational tools that can assess the system state and its operational margins more accurately and efficiently than current approaches. Specifically, these new techniques need to go beyond linearized methods of analysis and ensure that the power system is  safe even in the presence of large disturbances and uncertainty. In this paper, we study the energy function approach and argue that this classic tool can be generalized to help solve a number of existing and emerging problems in operating modern power systems.

Introduced more than 50 years ago by Aylett \cite{aylett1958energy}, energy functions have been used to develop direct methods of the power systems stability and analysis. Essentially, energy functions are physically motivated Lyapunov functions for certain simplified models of the power system dynamics. In the 70s and 80s, many results were derived for stability analysis of power systems,first via approximate Kron-reduction techniques and then using more accurate structure-preserving energy functions \cite{BergenHill1981}.  Developments of this period are summarized in \cite{varaiya1985direct,pai1989energy}. 
Since then, much of the work has focussed on developing effective heuristics to produce non-conservative estimates of the region of power system stability (see \cite{chiang2011direct} and references therein).
Interest in these approaches declined in the 1990s and 2000s, partially due to the difficulty of generalizing these approaches to more detailed power system models, but also due to advances in computing technologies that made stability analysis based on real-time dynamic simulations feasible.

Our results are also linked to more recent lines of research related to (a) necessary and sufficient conditions guaranteeing existence and uniqueness power flow solutions \cite{ilic1992network}\cite{araposthatis1981analysis}\cite{lesieutre1999existence}\cite{dorfler2014synchronization}\cite{dorfler2012synchronization}\cite{dorfler2013synchronization}; (b) analysis of the fixed point characterizations of the power flow equations and analytical approximation schemes for power flow solutions \cite{bolognani2014existence}; (c) semidefinite relaxations of optimal power flow \cite{low2014convex,low2014convexb} and various extensions (see for example\cite{louca2014nondegeneracy,madani2014promises}). We argue in this manuscript that energy functions can be viewed as a unifying tool for several operational problems in power systems. This is based on the fact that for lossless systems (resistance of transmission lines are ignored) stationary points of energy functions are solutions of the Power Flow (PF) equations. Hence, many traditional operational problems dependent on the notion of the (static) PF equations, such as the PF analysis itself, the Critical Loadability Problem (CLP) and Optimum Power Flow (OPF) problem,can be studied using the energy function methodology.  This viewpoint, first presented in \cite{12BV,bent2013synchronization}, is extended here both in terms of theory and applications.
Moreover, we show that this viewpoint extends to the more general cases where the energy-like functions are not true Lyapunov functions, but simply functions whose stationary points correspond to PF solutions. This extension will allow us to deal with lossy systems that have fixed resistance to reactance ratio for all transmission lines.

In order to obtain efficient formulations, we study a property of the energy functions that, to the best of our knowledge, has previously received little attention: \emph{convexity}. Since power systems are designed to operate at a stable equilibrium point, the energy (Lyapunov) functions characterizing their behavior are obviously convex at least in a small vicinity of the stable equilibrium. However, to the best of our knowledge, the exact domain of convexity has never been characterized. The classical papers by Bergen and Hill \cite{BergenHill1981} on the structure preserving model of power systems, extended in \cite{narasimhamurthi1984generalized,VanCutsem1985} to account for effects of the reactive power flows and variable voltages, form a solid base for our analysis.
We study various assumptions under which the energy function is {\bf\it provably convex}. It is well-known that when all nodes are the so-called \PV~ nodes (active power consumption/production and voltages are fixed/known), the energy function is convex \cite{bent2013synchronization,12BV} as long as all phase differences are smaller than $90$ degrees. In this paper, we generalize this result to the networks with \PQ~nodes (i.e. nodes where active and reactive power injections/consumptions are fixed while voltage magnitudes are variable).

Our main technical result is a description of the domain over which the energy function is convex. We provide a nonlinear but convex matrix inequality condition on the voltage phases and magnitudes that guarantees convexity of the energy function. Moreover, we show that this condition provides an \emph{exact} characterization of the domain of convexity of the energy function for the case of acyclic networks. For mesh networks, in general, the matrix inequality condition provides an inner approximation of the convexity domain. The condition for convexity is quite natural and captures the intuition that if voltage magnitudes and voltage phases at neighboring nodes are ``not too far'' from each other then the underlying energy function is convex.

The rest of this paper is organized as follows. In section \ref{sec:Intro}, we describe the mathematical background and power system models. In section \ref{sec:Energy}, we first analyze convexity of the energy function jointly in voltage magnitudes and phases, and then describe how the convexity can be used to design a provably convergent algorithm solving the PF equations (Section \ref{sec:Pflow}). In section \ref{sec:Related}, we discuss related work and describe the advantages of the energy function methodology. In section \ref{sec:Num}, we provide numerical illustrations of our approach on the IEEE 14 and 118 bus networks. Section \ref{sec:Conclusions} contains conclusions and a brief discussion of the future work. Finally,  Appendices \ref{subsec:Proof} and \ref{subsec:Lossy} are reserved for the proof of the main theorem and a generalization to a special class of lossy networks.


\section{Modeling Power Systems}\label{sec:Intro}

In order to guarantee the existence of an energy function, we will need to make the standard assumptions used in the structure preserving models \cite{BergenHill1981}, which are fairly well justified for the high-voltage transmission network:
\begin{itemize}
\item[(1)] All transmission lines are purely reactive (network is lossless), i.e. the network resistive power losses are considered small and thus ignored.
\item[(2)] Nodal active and reactive power injections/consumptions are constant.
\end{itemize}

\newcommand{\neb}{\sim}
\subsection{Notation}
The transmission network is modeled as a graph $\br{\Vs,\E}$ where $\Vs$ is the set of nodes and $\E$ is the set of edges. In power systems parlance, the nodes are called buses and the edges are called lines (transmission lines). We shall use these terms interchangeably in this manuscript. Nodes are denoted by indices $i=1,\ldots,n=|\Vs|$ and edges by ordered pairs of nodes $\br{i,j}$. We pick an arbitrary orientation for each edge, so that for an edge between $i$ and $j$, only one of $\br{i,j}$ and $\br{j,i}$ is in $\E$. If there is an edge between buses $i$ and $j$, we write $i \neb j,j\neb i$. Sometimes, it will be convenient to denote edges by a single index $k$. In this notation, $k_1$ and $k_2$ denote the ``from'' and ``to'' ends of the edge $k$. Define the edge-incidence matrix $A\in \R^{|\Vs|\times |\E|}$ as follows:
\begin{align}
A_{ik} =
\begin{cases}
1 \text{ if } k\in\E,i\in\Vs,k_1=i\nonumber\\
-1\text{ if } k\in\E,j\in\Vs,k_2=j\in\Vs\nonumber\\
 0 \text{ otherwise } \nonumber
       \end{cases}
       \end{align}

The set of vertices includes all load buses as well as buses representing generators (a single bus, or terminal and internal bus per generator). This allows us to develop a unifying framework that works for different load and generator models. Edges correspond to transmission lines. Each transmission line is characterized by its admittance  $G_{ij}-\ic B_{ij}$ ($\ic=\sqrt{-1}$). However, since we assume that transmission lines are lossless, we set $G_{ij}=0$ and thus each line is characterized by a single real and positive parameter $B_{ij}$ (the negated line susceptance).  

Let $V_i, \theta_i, P_i$ and $Q_i$ denote voltage (magnitude), phase, active and reactive injection at the bus $i$ respectively. Let $\rho_i=\logb{V_i}$. Let $\theta_{ij}=\theta_i-\theta_j,\rho_{ij}=\rho_i-\rho_j$. 
Buses are of three types:
\begin{itemize}
\item\underline{\PV buses} where active power injection and voltage are fixed, while voltage phase and reactive power adjust as conditions (e.g. power flows) are variables. The set of \PV buses is denoted by $\G$.
\item\underline{\PQ buses} where active and reactive power injections are fixed, while voltage phase and magnitude are variables. The set of \PQ buses is denoted by $\Lo$.
\item\underline{Slack bus}, a reference bus at which the voltage magnitude and phase are fixed, and the active and reactive power injections are free variables. The slack bus is denoted by $\Sb$.
\end{itemize}
$\theta$ and $V$ (with no subscript) denote the vectors of phases at all buses except the slack bus and voltage magnitudes at all $\PQ$ buses. $\rho,V$ are indexed by $\{i:i\in \Lo\}$. For every bus $i$ define $
B_i=\sum_{j\sim i}B_{ij}$. 

Let $S$ be an ordered set of indices. We use the notation
\[M=\brs{\begin{pmatrix}a & b\\ c & d\end{pmatrix}}_{ij}^{S}\]
where $i,j\in S$, to denote an $|S| \times |S|$ matrix whose rows and columns are indexed by the entries of $S$, with $M_{ii}=a,M_{ij}=b,M_{ji}=c,M_{jj}=d$ and all other entries equal to $0$. Similarly, we use the notation $M=[a]_{ii}^{S}$ to denote an $|S|\times |S|$ matrix with the $(i,i)-th$ entry equal to $a$ and all other entries equal to $0$.

Given a vector $x$ with indices $i\in S$, denote by $\mathrm{diag}\br{x}$ the $|S|\times|S|$ matrix with diagonal entries $x_i$. 
Given a set $\C\subset\R^n$, $\Int\br{\C}$ denotes the interior of the set.

\subsection{Background}

In this section, we introduce two power systems concepts that are critical for what follows: PF equations and energy functions.

\subsubsection{Power Flow Equations}

PF Equations model the flow of power over the power system network. They are a set of coupled nonlinear equations that follow from Kirchoff's laws applied to the AC power network. Circuit elements in the standard power systems models are all linear, if one ignores discrete elements like phase shifters and tap-changing transformers. Even though Ohm's law and Kirchhoff's law are linear in voltages and currents, power is their product and hence quadratic.   Since we are interested in power, e.g., to determine power generation, specify demand, minimize line loss, PF equations are nonlinear.

PF equations assume that the network is balanced, that is, the net sum of powers consumptions, injections and power dissipated is zero. This relies on the assumption that at the time-scale where the PF are solved (every 5 minutes or so), the system is in a quasi-steady state, i.e. the dynamic disturbances have been resolved through actions of the automatic control schemes represented by the automatic voltage regulators, power system stabilizers and primary and secondary frequency control systems.

Since our main concern here is to study properties of the PF equations, we simply state them without derivation. At each node in the power system, there are four variables: voltage magnitude ($V_i$), voltage phase ($\theta_i$), active power injection ($P_i$), reactive power injection ($Q_i$). When solving the PF equations, at each node, two of these variables are fixed and the other two are free. This leads to $2n$ equations in $2n$ variables, where $n$ is the number of buses.
Typically, power systems have two types of nodes: \PV~nodes and \PQ~nodes.
At the \PV~buses, the active power injection $P_i$ and the voltage magnitude $V_i=\expb{\rho_i}$ are fixed. The fixed voltage magnitude is denoted $\Vset_i$, to distinguish it from variable voltages at other nodes. Generators are typically modeled as \PV~buses, since they have voltage regulators that inject reactive power to maintain voltages at fixed values. 
We can choose $\Vset_i=1$ for each $i\in\G$, since any non-unit voltage set-point can be absorbed into the susceptances $B_{ij}$ in the active and reactive power flow equations.

At the \PQ~buses, the active power injection $P_i$ and the reactive power injection $Q_i$ are held fixed. Loads are typically \PQ~nodes, where $P_i$ and $Q_i$ represent the active and reactive power demands (assumed constant). 

Usually, the variables $Q_i$ at the \PV~buses are not considered explicitly. These variables can be obtained from the reactive power balance equations by plugging in there the values of $\br{\rho,\theta}$. Finally, we have a special bus (slack bus) which absorbs all the power imbalances in the network.
The injections $P_{\Sb},Q_{\Sb}$ at this bus are free variables whose values are adjusted to guarantee the power balance and $\theta_{\Sb}=\rho_{\Sb}=0$.

The power balance equations at all the buses,   a set of coupled nonlinear equations in $\br{\rho,\theta}$, together constitute the PF equations:
\begin{subequations}
\begin{align}
 P_i & =\sum_{j \neb i} B_{ij}\expb{\rho_i+\rho_j}\sin\br{\theta_{ij}} ,i\in\G\cup\Lo\label{eq:pva}\\
\rho_i & =0 ,i\in\G\label{eq:pvb}\\
Q_i &=\sum_{j \neb i} B_{ij}\br{\expb{2\rho_i}-\expb{\rho_i+\rho_j}\cos\br{\theta_{ij}}} ,i\in\Lo\label{eq:pqb}\\
\theta_{\Sb} &=0 ,\rho_{\Sb}=0.\label{eq:sab}
\end{align}\label{eq:pf}
\end{subequations}

\subsubsection{Energy Functions}

Energy functions were first introduced in the context of first integral analysis of the power system swing dynamic equations \cite{BergenHill1981}\cite{VanCutsem1985}. The energy function consists of two terms -- kinetic and potential.  It is also a Lyapunov function for the swing dynamics.
In this paper, we are not explicitly concerned with the dynamics, so the most useful property of the energy functions is that the stationary points of the energy function potential part map to solutions of the PF equations.  The PF equations \eqref{eq:pf} can be  re-written in the following variational form
\begin{subequations}
\begin{align}
\text{\eqref{eq:pva} : }   0
& =\frac{\partial E}{\partial \theta_i} \quad \forall i\in\Lo\cup\G\label{eq:pvvara}\\
\text{ \eqref{eq:pqb} : }   & 0= \frac{\partial E}{\partial \rho_i}\quad \forall i\in\Lo\label{eq:pqvarb}\\
& \rho_i=0\quad \forall i\in\G, \rho_{\Sb}=\theta_{\Sb}=0
\end{align}\label{eq:pfvar}
\end{subequations}
where $E$ is called the energy function:
\begin{align}
E\br{\rho,\theta} & = -\sum_{i\in\G\cup\Lo}P_i \theta_i -\sum_{i\in\Lo} Q_i\rho_i \nonumber\\
& \, + \sum_{\br{i,j}\in\E} B_{ij}\br{\frac{\expb{2\rho_i}+\expb{2\rho_j}}{2}}\nonumber\\
&\, -\sum_{\br{i,j}\in\E} B_{ij}\expb{\rho_i+\rho_j}\cos\br{\theta_{ij}}\label{eq:Energyrho}
\end{align}
Note that $\rho_i$ is a variable for $i \in \Lo$ and $\rho_i=0$ for $i \in \G$, $\rho_{\Sb}=0,\theta_{\Sb}=0$.

\section{Energy Function and Convexity}\label{sec:Energy}

In this section, we study the convexity of the energy function \eqref{eq:Energyrho}. We first discuss why convexity is an interesting and useful property, and then describe our results characterizing the domain of convexity of the energy function \eqref{eq:Energyrho}. Before we proceed, let us define the domain of convexity precisely.
\begin{definition}
The domain of convexity of the energy function is defined as:
\[\D=\left\{\br{\rho,\theta}:\nabla^2 E\br{\alpha \rho,\alpha \theta}\succeq 0 \quad \forall \alpha\in[0,1]\right\}\]
\end{definition}

\begin{remark}
In general, there could be several (disconnected) regions over which the energy function is convex \cite{araposthatis1981analysis}. We require that for any $\br{\rho,\theta}\in\D$, the energy function is convex at every point on the line segment connecting the origin to $\br{\rho,\theta}$. This ensures that we pick the region of convexity that includes the origin.
\end{remark}

\subsection{Why is Convexity Interesting?}

The major contribution of this paper is the characterization of the convexity domain of the energy function \eqref{eq:Energyrho}. We now discuss why this characterization is important and interesting. The principal reason is that solutions of the PF equations are stationary points of the energy function \eqref{eq:pfvar}. We will show, in section \ref{sec:Pflow}, that one can construct a convex optimization problem that will either find a solution of the PF equations within the convexity domain, or otherwise certify that there exist no solutions within the convexity domain. In this Section, we justify restricting ourselves to PF solutions contained in $\D$, as the solutions in $\D$ show some additional desirable properties.

\subsubsection{Asymptotic Stability}

The dynamics of the power system, under certain reasonable assumptions, reduce to the so-called ``swing dynamics'', see e.g. \cite{varaiya1985direct}.  In this model, each generator is modeled as a constant voltage source behind a transient reactance.  For each generator, there is an ``internal'' node (modeled as a constant voltage source subjected to the second order angular dynamics) and a ``terminal'' node (a node with $0$ active and reactive power injections). We denote the set of internal generator nodes by $\G_0$. The swing equation at an internal node $i\in\G_0$ is
\begin{align}
M_i\ddot{\theta_i}+D_i\dot{\theta_i}=-\frac{\partial E\br{\rho,\theta}}{\partial \theta_i} \label{eq:pvdyn}
\end{align}
At all other nodes $i\in\G\cup\Lo$, we have the \emph{algebraic} active and reactive power balance equations \eqref{eq:pf}. The dynamic state variables are $\{\dot{\theta_i}:i\in\G_0\}\cup\{\theta_i:i\in\G_0\}$.
\begin{thm}
Let $\br{\rho,\theta}$ be an equilibrium of the swing dynamics such that $\br{\rho,\theta}\in\D$ and satisfy all the algebraic equations \eqref{eq:pf}. Then, $\br{\rho,\theta}$ is asymptotically stable, that is, linearizing \eqref{eq:pvdyn} around $\br{\rho,\theta}$ results in a stable linear system.
\end{thm}
\begin{proof}
This follows from Lemma 5.1 in \cite{tsolas1985structure}, which is a stronger result.
\end{proof}

\subsubsection{Existence of Solutions}\label{sec:SolExist}

The following result is contained in our recent manuscript \cite{DjMallada}:
\begin{thm}
Suppose that the network is a tree. Then, the PF equations \eqref{eq:pf} have a solution if and only if  a solution is contained within $\D$.
\end{thm}
Thus, at least for trees, it is sufficient to look for solutions of the PF equations contained in $\D$. Preliminary numerical evidence (section \ref{sec:Num}) suggests that this observation also extends to mesh networks. Further, as we shall illustrate numerically in Section \ref{sec:Num}, the solutions within $\D$ are the ``desirable'' solutions (they have small phase differences between neighboring buses and sufficiently large voltage magnitudes at all \PQ~buses). 


\subsection{Characterization of Domain of Convexity of the Energy Function}

We now derive our main results on convexity of the energy functions. Since stationary points of the energy functions correspond to solutions of the PF equations, which are known to have multiple isolated solutions, the energy function is not globally convex. Thus,  we focus on searching for restricted domains within which the energy functions are convex. In particular, we will characterize a subset of $\D$ such that phase differences between neighbors are smaller than $90\deg$, since this is typical for normal power systems operations. 
It was shown in \citep{bent2013synchronization} that this restriction is actually sufficient for convexity in the case when all buses are \PV. However,  when the network also contains \PQ~buses additional conditions are required.

\begin{thm}\label{thm:EConvex}
The energy function $E\br{\rho,\theta}$ is jointly convex in $\br{\rho,\theta}$ over the convex domain $\C\subseteq\D$ given by
\begin{align}
&|\theta_i-\theta_j|\leq \frac{\pi}{2}\quad\forall (i,j)\in\E\label{eq:ConvC1}\\
&\sum_{i\in\Lo} \left[2B_i-\sum_{j\in\G\cup\{\Sb\}}B_{ij}\frac{\expb{\rho_{ji}}}{\cos\br{\theta_{ji}}}\right]_{ii}^{\Lo}\nonumber \\
& \,-\sum_{(i,j)\in\E,i,j\in\Lo}\frac{B_{ij}}{\cos\br{\theta_{ij}}} \brs{\begin{pmatrix}\expb{\rho_{ji}} & 1 \\ 1 & \expb{\rho_{ij}}\end{pmatrix}}_{ij}^{\Lo}\succeq 0\label{eq:ConvC2}
\end{align}
Further, suppose that the network has a tree topology. Then, the domain of function convexity, $\D$, equals $\C$. Also, the energy function is \emph{strictly convex} on $\Int\br{\C}$.
\end{thm}
\begin{proof}
See appendix section \ref{sec:App}.
\end{proof}A number of Remarks are in order.
\begin{remark}{\it Simple special cases}\\
When no \PQ~nodes are connected to each other, the second term of the semi-definite inequality \eqref{eq:ConvC2} is an empty sum. Thus, the semi-definite inequality is equivalent to
\[2B_i=2\br{\sum_{j\sim i} B_{ij}} \geq \sum_{j\in\G\cup\{\Sb\},j\sim i}B_{ij}\frac{\Vset_j\expb{-\rho_i}}{\cos\br{\theta_i-\theta_j}}\]
Since $B_i=\sum_{j\sim i} B_{ij}$, it suffices that $\Vset_i\expb{-\rho_j}\leq 2\cos\br{\theta_i-\theta_j}$ for this to be true, which is satisfied if $V_j\geq .7 \Vset_i$, $|\theta_i-\theta_j|\leq \frac{\pi}{4}$. This is a fairly reasonable restriction in practice. A similar condition was discussed in \cite{tsolas1985structure}. However, it was imposed only at the equilibrium point where the energy function reaches its minimum,  as a sufficient condition to guarantee (small deviation) asymptotic stability. The asymptotic stability follows from our results as well, which guarantees that the equilibrium point lies within the domain of convexity of the energy function and it is hence asymptotically stable. 

Also, if there are no \PQ~nodes, the matrix inequality condition in \eqref{eq:ConvC2} is empty and hence the convexity condition reduces to requiring that all phases differences (over existing lines) are smaller than $\frac{\pi}{2}$ (see \cite{bent2013synchronization}).

When \PQ-\PV~connections are present, the semidefinite inequality does not simplify. However, in section \ref{sec:Num}, we show numerically that the conditions imposed by the convexity domain are not very restrictive and most practical power flow solutions lie within the convexity domain.
\end{remark}

\begin{remark}{\it Conservatism} \\
For the tree networks, we have shown that the convexity conditions are necessary and sufficient (assuming that we are interested only in the domain with phase differences over all the lines smaller than $\frac{\pi}{2}$). For meshed networks, these conditions are sufficient but, possibly, not necessary. This is due to the fact that our convexity analysis has treated the edge variables $\theta_{ij}$ as  independent for each edge $(i,j)\in\E$ while in reality the variables are coupled since $\theta_{ij}=\theta_i-\theta_j$. Analysis of the gap is left for future work. We present some numerical results indicating that the gap is small in section \ref{sec:Num}.
\end{remark}

\subsection{Solving Power Flow Equations via Convex Optimization}\label{sec:Pflow}

Although energy function was developed initially as a tool to assess dynamic stability, the convexity results also yields consequences useful in the context of  a number of other important power system applications. We choose to emphasize here one of these applications: Solution of the PF equations, leaving the discussion of other potential applications for future work.

As we have already seen, the energy function technique provides a variational characterization of the PF equations \eqref{eq:pfvar}. Therefore, finding solution of the PF equations is equivalent to finding a stationary point of the energy function.

\begin{corollary}
Let $S\subset \Int\br{\C}$. Then, the power flow equations \eqref{eq:pf} have a solution $\br{\opt{\rho},\opt{\theta}}\in S$ if and only if the following optimization problem has its (unique) optimal solution in $S$:
\begin{align}
\min_{\br{\rho,\theta}\in \C,\theta_{\Sb}=0,\rho_{\Sb}=0} E\br{\rho,\theta}
\label{eq:optE}
\end{align}
\end{corollary}
\begin{proof}
If $\br{\opt{\rho},\opt{\theta}}\in S$ solves \eqref{eq:pf}, it satisfies $\nabla_{\rho}E\br{\opt{\rho},\opt{\theta}}=0,\nabla_{\theta}E\br{\opt{\rho},\opt{\theta}}=0$. Thus, it also satisfies the KKT conditions for the optimization problem \eqref{eq:optE}, forcing the respective Lagrange multipliers to be $0$. Hence, the solution is optimal for \eqref{eq:optE}.

Conversely, suppose that $\br{\opt{\rho},\opt{\theta}}\in  S\subset\mathrm{int}\br{\C}$ solves \eqref{eq:optE}. Then, by complementary slackness, the Lagrange multipliers are $0$ and hence the KKT conditions reduce to $\nabla_{\rho}E\br{\rho,\theta}=0,\nabla_{\theta}E\br{\rho,\theta}=0$. Thus, $\br{\opt{\rho},\opt{\theta}}\in\C$ also solves \eqref{eq:pf}.

Further, by strict convexity of the energy function $E\br{\rho,\theta}$ on $\Int\br{\C}$, there can be at most one stationary point of $E$ in $\Int\br{\C}$, and hence, the power flow solution, if it exists, is unique.
\end{proof}

\section{Related Work}

In this section, we discuss and contrast our work with related previous work.
\subsection{Principal Singular Surfaces}
The work closest in spirit to our approach is presented in the series of papers \cite{TavoraA,TavoraB,TavoraC}. The authors of \cite{TavoraA,TavoraB,TavoraC} consider lossless networks where all buses are \PV~buses. They define Singular surfaces, which are closed surfaces where $\detb{\nabla^2 E}=0$. The Principal Singular Surface (PSS) is the unique surface that encloses the origin $\theta=0$ and the Principal Region is the region enclosed by it. The principal region coincides with our definition of $\D$ when all nodes are \PV~nodes, apart from the additional requirement that the phase differences are smaller than $\frac{\pi}{2}$. Our results characterizing $\C$ can be seen as constructing convex inner approximations to the principal region.

This analysis was extended in \cite{araposthatis1981analysis}, where a number of counterintuitive properties of the set of solutions of the power flow equations were shown , even in the simple case when all buses are \PV~buses. Based on these results, it was conjectured in this paper that the principal region is nonconvex in general. Further, there are examples of networks for which the only solutions to the power flow equations are unstable solutions. However, it was conjectured that if there is a stable solution, there must be one in the principal region. 

Extending this analysis and verifying the conjectures for our general setting with \PQ~nodes is a direction of future work.

\subsection{Convex Relaxations of OPF}

The authors of the series of recent papers \cite{lavaei2013geometry,zhang2013geometry,lavaei2013geometry,low2014convex,low2014convexb} have studied \emph{convex relaxations} of Optimal Power Flow (OPF) over tree networks. 
Various conditions were uncovered under which the convex relaxations of the OPF problems are exact, in the sense that an optimal solution of the original OPF problem can be recovered from the solution of the convex relaxation. In this work, we primarily deal with the PF equations, not OPF. However, for lossless networks our approach applies to arbitrary topologies and can potentially be extended to the OPF formulation, following the logic first time sketched in \cite{bent2013synchronization}:
\begin{align*}
\max_{P,Q \in T} \min_{\br{\rho,\theta}\in S} & -\sum_{i\in\G\cup\Lo}P_i \theta_i -\sum_{i\in\Lo} Q_i\rho_i \nonumber \\
&  +\sum_{\br{i,j}\in\E} B_{ij}\br{\frac{\expb{2\rho_i}+\expb{2\rho_j}}{2}}\nonumber\\
& -\sum_{\br{i,j}\in\E} B_{ij}\expb{\rho_i+\rho_j}\cos\br{\theta_{ij}}\nonumber \\
& \quad -\lambda \costt\br{P,Q}
\end{align*}
where $\cost\br{P,Q}$ is a convex cost on the injections and $T,S$ are convex operational constraint sets on the injections and voltages, respectively and $\lambda>0$ is a positive scaling factor. This is a convex-concave saddle point problem, as long as $S\subset \C$ and can be solved efficiently using ideas from duality \cite{boyd2009convex}. If the constraints in $S$ are not binding at the optimal solution, the inner minimization effectively enforces the PF equations \eqref{eq:pf}, since the minimization corresponds to setting the gradient of the energy function to $0$. The outer maximization then optimizes a cumulative cost consisting of the negative injection cost with an extra term added -- energy function value at the optimum over line flows (corresponding to solution of the PF). We can choose $\lambda$ large enough so that the energy function term can be neglected for the outer maximization, so that we are effectively solving the OPF problem:
\begin{align*}
\min_{P,Q,\rho,\theta} & \costt\br{P,Q}\\
\text{Subject to } & \text{\eqref{eq:pf}}, \br{P,Q}\in T,\br{\rho,\theta} \in S
\end{align*}
Further developments on this line of work, along with conditions when our minimax formulation is equivalent to the OPF problem above, will be pursued in subsequent work.

\subsection{Conditions on Existence of solutions to PF equations}

Several papers have studied conditions for existence of solutions to the Power Flow Equations \cite{bolognani2014existence}\cite{lesieutre1999existence}\cite{wu1982steady}. In \cite{molzahn2012sufficient}, the authors propose a sufficient condition for the insolvability of power flow equations based on a convex relaxation.
However, our approach differs from these in the following important ways:

\begin{itemize}
\item To solve the PF equations in $\C$, we provide necessary and sufficient conditions, that is, our approach finds a solution in $\C$ if and only if there exists a solution in $\C$.

\item Our approach is algorithmic, that is, we provide an algorithm (based on convex optimization) that is guaranteed to find the solution efficiently (in polynomial time).

\item If there are additional operational constraints $\br{\rho,\theta}\in S$, with $S\subset \C$, e.g. correspondent to line protection \cite{singh2001direct} and/or line flow limits, we can additionally answer the question of whether there exists a PF solution with $\br{\rho,\theta}\in S$. This is an important contribution, since most of the time system operators are interested in finding power flow solutions that additionally satisfy operational constraints. \end{itemize}

In \cite{bolognani2014existence}, the authors also propose an algorithm based on a contraction mapping consideration. However, the algorithm only works in a small ball around the origin in the $\br{P,Q}$ space. Our results are stated in terms of a nonlinear convex constraint in $\br{\rho,\theta}$ space. Precisely understanding the set of $\br{P,Q}$ for which the solution $\br{\rho,\theta}\in\C$, and conversely the set of $\br{\rho,\theta}$ for which $\br{P,Q}$ lies in a certain ball, is still an open problem, even for the special case where all buses are \PV~buses. This setting was studied and the results were exteded and connected with results on synchronization in coupled oscillators in a series of recent papers \cite{dorfler2012synchronization}\cite{dorfler2013synchronization}\cite{dorfler2014synchronization}. The authors provide distinct sufficient and necessary conditions on the injections for the existence of power flow solutions with phase differences satisfying certain bounds. The authors also note that obtaining non-conservative sufficient conditions is still an open problem, and show numerically that a necessary condition is ``almost'' sufficient for most power systems. Extending this analysis to the setting with \PQ~buses and obtaining non-conservative sufficient conditions on the injections that guarantee existence of power flow solutions is an interesting direction for future work.


\section{Numerical Illustrations}\label{sec:Num}

This Section presents numerical experiments illustrating theoretical results resented above.

\subsection{Existence of Solutions: 2-bus Network}

Here we analyze the toy case of a two bus network with one generator and one load. The generator also serves as the slack bus. Thus, the only variables are the voltage phase $\theta$ and voltage magnitude $\expb{\rho}$ at the load bus, and the parameters are the active power demand $P$ and the reactive power demand $Q$ at the load bus. The results, plotted in the figure \ref{fig:Convexity} show the power flow solutions shown on the heat (iso-line) map of the energy function. There are two solutions for small loads. The convexity domain is the region inside the dashed green arc. The solution plotted as a green dot is a dynamically stable solution in the interior of the convexity domain. The solution plotted as a red dot is a dynamically unstable solution outside the convexity domain. As the loads are increased, the two solutions move closer to each other. At a certain critical value, the two solutions hit the boundary of the convexity domain and then vanish, thus showing that, beyond this critical loading level, there are no longer any solutions to the PF equations.

This simple experiment supports our conjecture that in general, if there is a solution to the PF equations, there must be one found within the convexity domain. As stated earlier, this result was proven for networks with tree topology \cite{DjMallada}. In the special case when all nodes are the \PV~buses, this result was conjectured (without proof) in \cite{araposthatis1981analysis}, with the caveat that one is only looking for \emph{dynamically stable} solutions. We plan to investigate this conjecture for the general case with \PQ~buses in further work.

\begin{figure}[htb]
        \centering
        \begin{subfigure}[b]{0.32\columnwidth}
                \includegraphics[width=\textwidth]{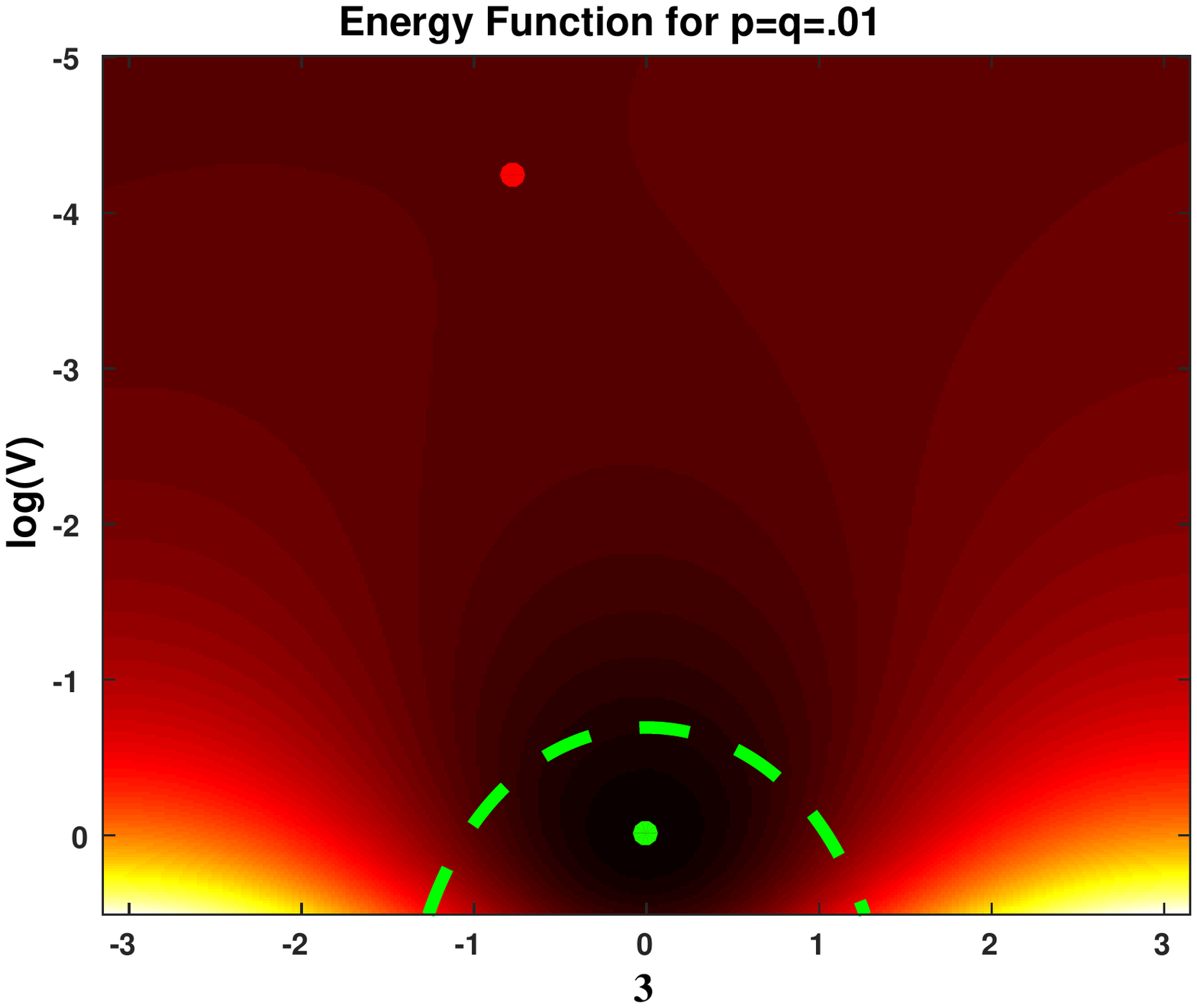}
                \caption{$P=Q=.01$}~\label{fig:TwoA}
        \end{subfigure}
        \begin{subfigure}[b]{0.32\columnwidth}
                \includegraphics[width=\textwidth]{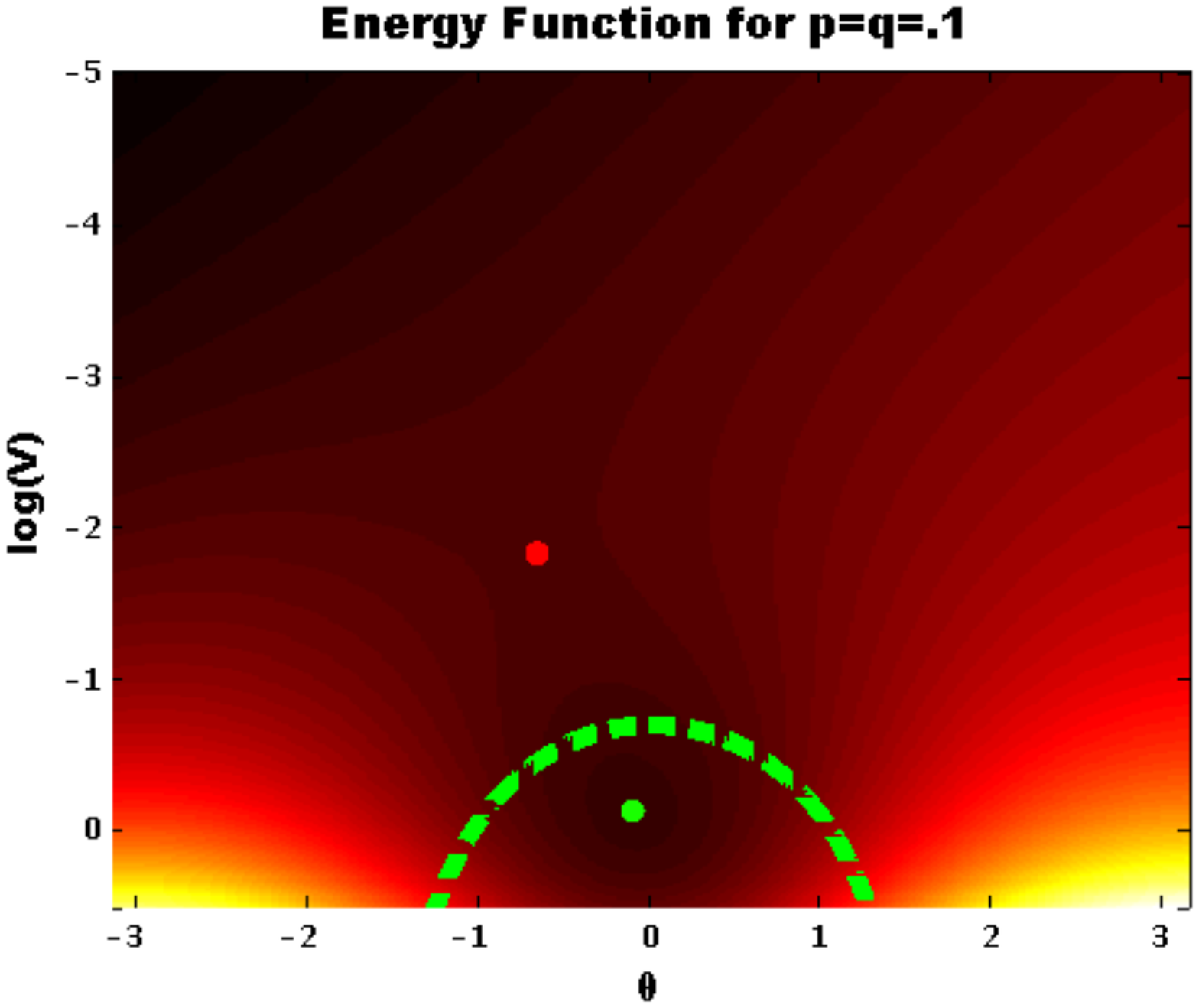}
                \caption{$P=Q=.1$}~\label{fig:TwoB}
        \end{subfigure}\\
        \begin{subfigure}[b]{0.32\columnwidth}
                \includegraphics[width=\textwidth]{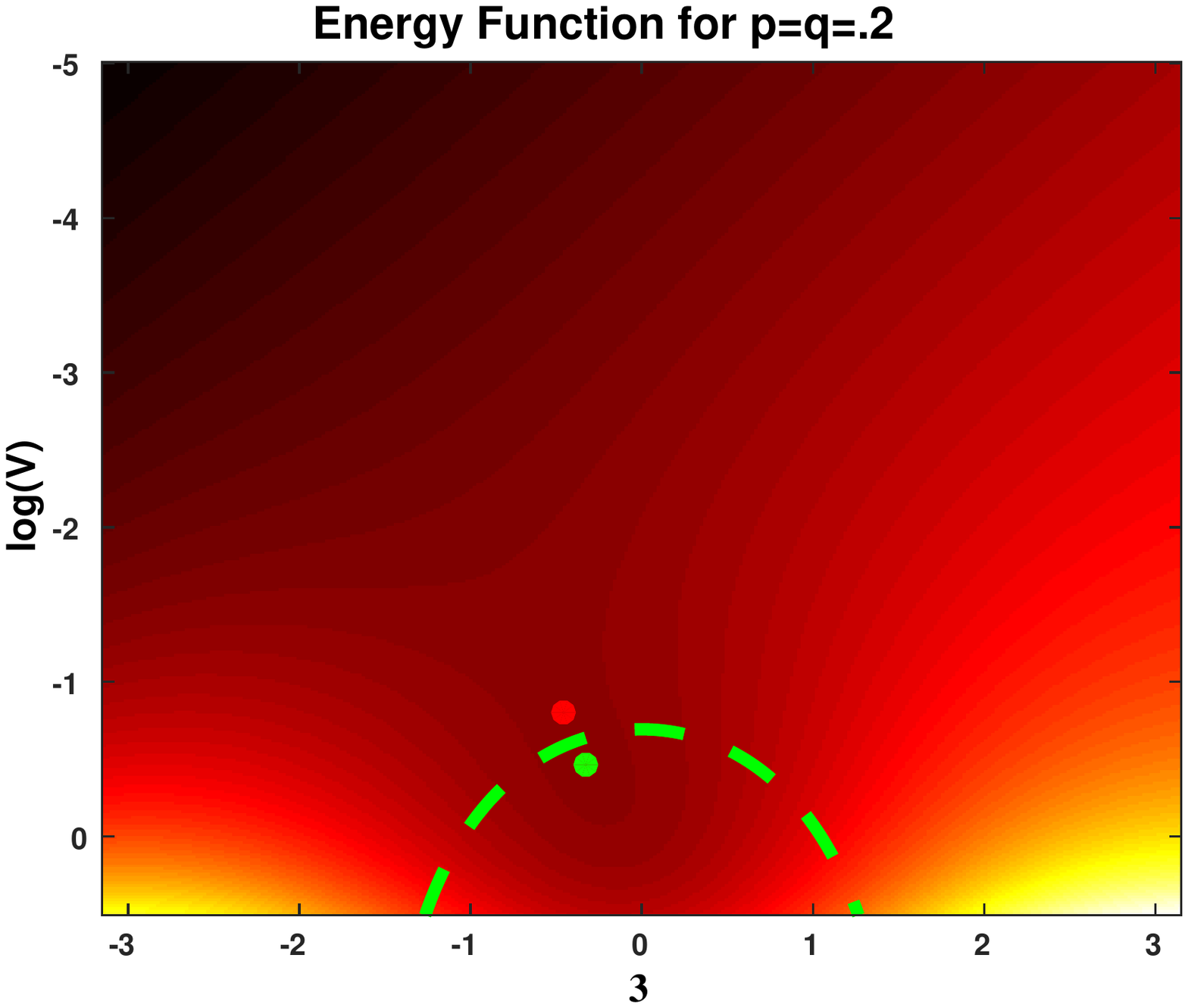}
                \caption{$P=Q=.2$}~\label{fig:TwoC}
        \end{subfigure}
         \begin{subfigure}[b]{0.32\columnwidth}
                \includegraphics[width=\textwidth]{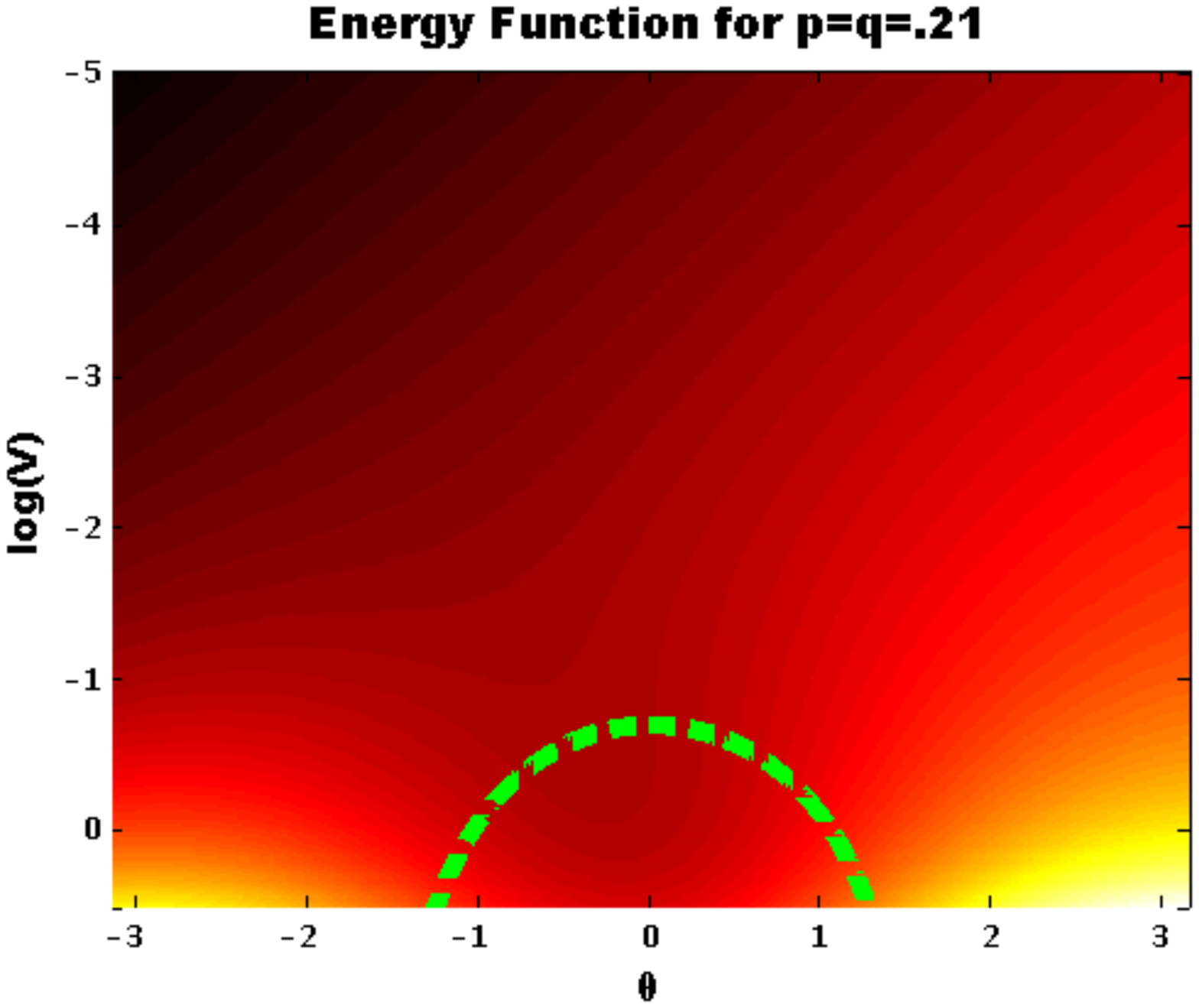}
                \caption{$P=Q=.25$}~\label{fig:TwoD}
        \end{subfigure}

        \caption{Solutions of the PF equations (red and green dotes) shown within the heat (iso-line) map of the energy function. The green dot is a dynamically stable solution inside the convexity domain, and the red dot is a dynamically unstable solution outside the convexity domain. Two distinct solutions are sufficiently far from each other at the initial (moderate) load. The two solutions get closer to each as the load increases. The solutions approach the boundary of the region of convexity from opposite sides and disappear once they hit the boundary.}\label{fig:Convexity}
\end{figure}
\subsection{Theoretical versus Numerical Domain of Convexity}

In this section, we compare numerically these regions of the energy function and the reduced energy function convexity. In order to be able to visualize these regions, we restrict ourselves to a 3 bus toy system with one PV bus and two \PQ buses. Bus 1 is the PV bus (also the slack bus): its voltage magnitude $V_1$ is fixed to 1 pu and its phase is selected as a reference, $\theta_1=0$. Buses 2 and 3 are \PQ buses with variable voltage magnitude $V_2=\expb{\rho_2},V_3=\expb{\rho_3}$ and voltage phase $\theta_2,\theta_3$. $V_2,V_3$ are determined by solving the reactive PF equations
\begin{align*}
& Q_2=B_{12}\br{V_2^2-V_2\cos\br{\theta_2}}+B_{23}\br{V_2^2-V_2V_3\cos\br{\theta_2-\theta_3}} \\
& Q_3=B_{13}\br{V_3^2-V_3\cos\br{\theta_3}}+B_{23}\br{V_3^2-V_2V_3\cos\br{\theta_2-\theta_3}}
\end{align*}
Active power injections are chosen as $P_1=3.7,P_2=-1.7,P_3=-2$  (these numbers are in the normalized units) and the reactive power demands are $Q_2=-1.05,Q_3=-1.24$.  The susceptances on the lines are $B_{12}=26.88,B_{13}=26.88$. We consider two cases for $B_{23}$ that lead to different network topologies: In the first case buses $2,3$ are not connected so $B_{23}=0$ and in the second case we choose $B_{23}=16.67$.

Given $\theta$, we can solve the reactive power flow equations \eqref{eq:pqb} for $\rho$: We denote the solution by $\rho\br{\theta}$. 
We study the convexity domain of the reduced energy function $E\br{\rho\br{\theta},\theta}$ so that we can plot the results (since the reduced energy function is only a function of two variables, $\br{\theta_1,\theta_2}$).

For different values of $\theta_2,\theta_3\in\left[-\frac{\pi}{3},\frac{\pi}{3}\right]$, we compute the value of the energy function $E$ by solving the reactive PF equations for $\rho_2,\rho_3$ (using the Newton-Raphson iterative algorithm) and then plug  the results into $E\br{\rho,\theta}$. For any value of $\theta_2,\theta_3$, we also estimate the Hessian $\nabla^2_{\theta}E\br{\rho(\theta),\theta}$ by means of numerical differentiation. This gives us the exact region over which the reduced energy function, $E\br{\rho(\theta),\theta}$, is a convex function of $\theta$. We compare this to the domain $\C$ of the energy function (and thus reduced energy function) convexity predicted by our theorems. Since the theorems only predict joint convexity of $E\br{\rho,\theta}$, we deduce reduced convexity by finding the set of $\theta$ such that $\br{\rho\br{\theta},\theta}\in \C$.
\begin{figure}[htb]
        \centering
        \begin{subfigure}[b]{0.24\columnwidth}
                \includegraphics[width=\textwidth]{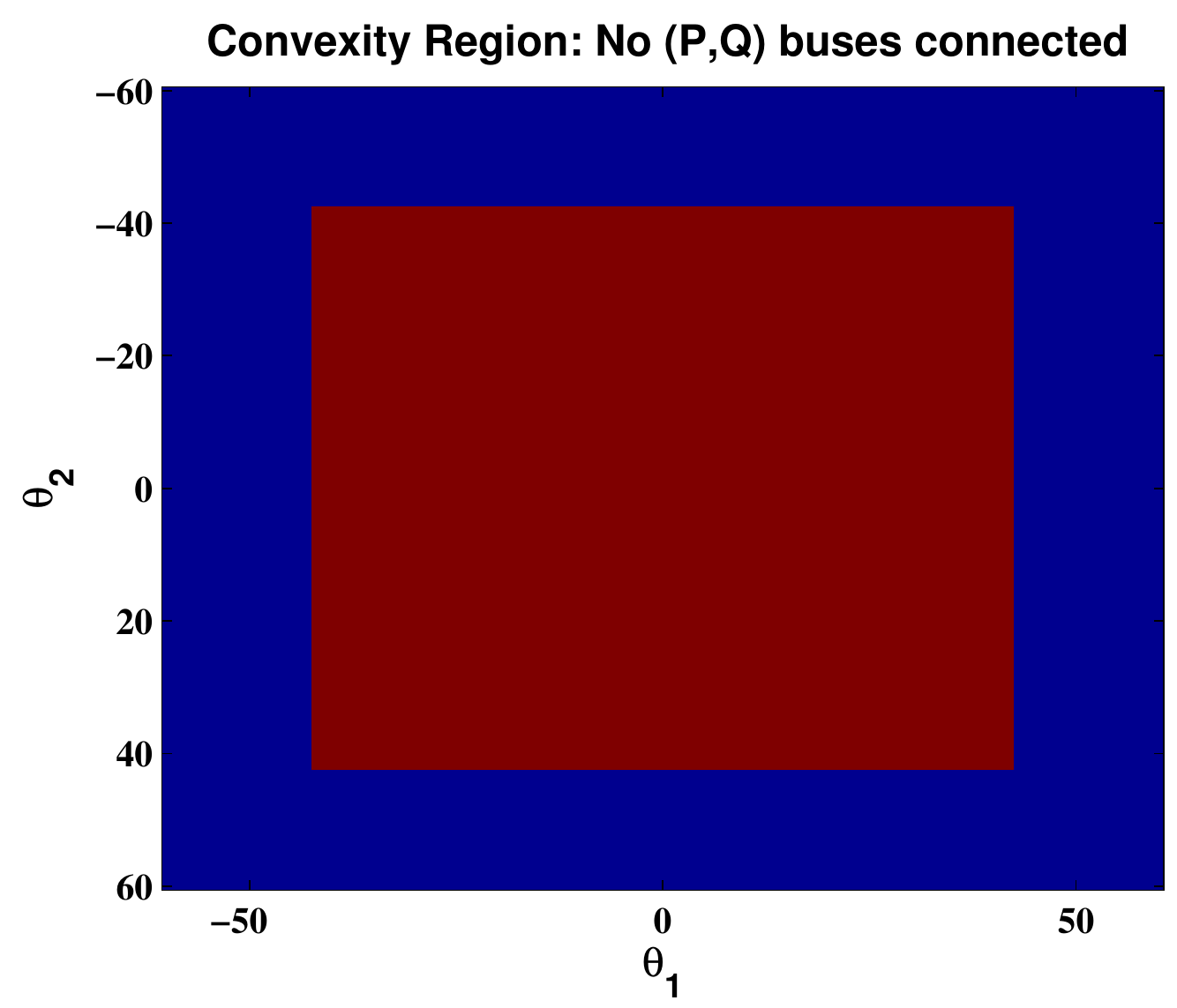}
                \caption{No (P,Q) connected}~\label{fig:ConvexityNoLoop}
        \end{subfigure}
        \begin{subfigure}[b]{0.24\columnwidth}
                \includegraphics[width=\textwidth]{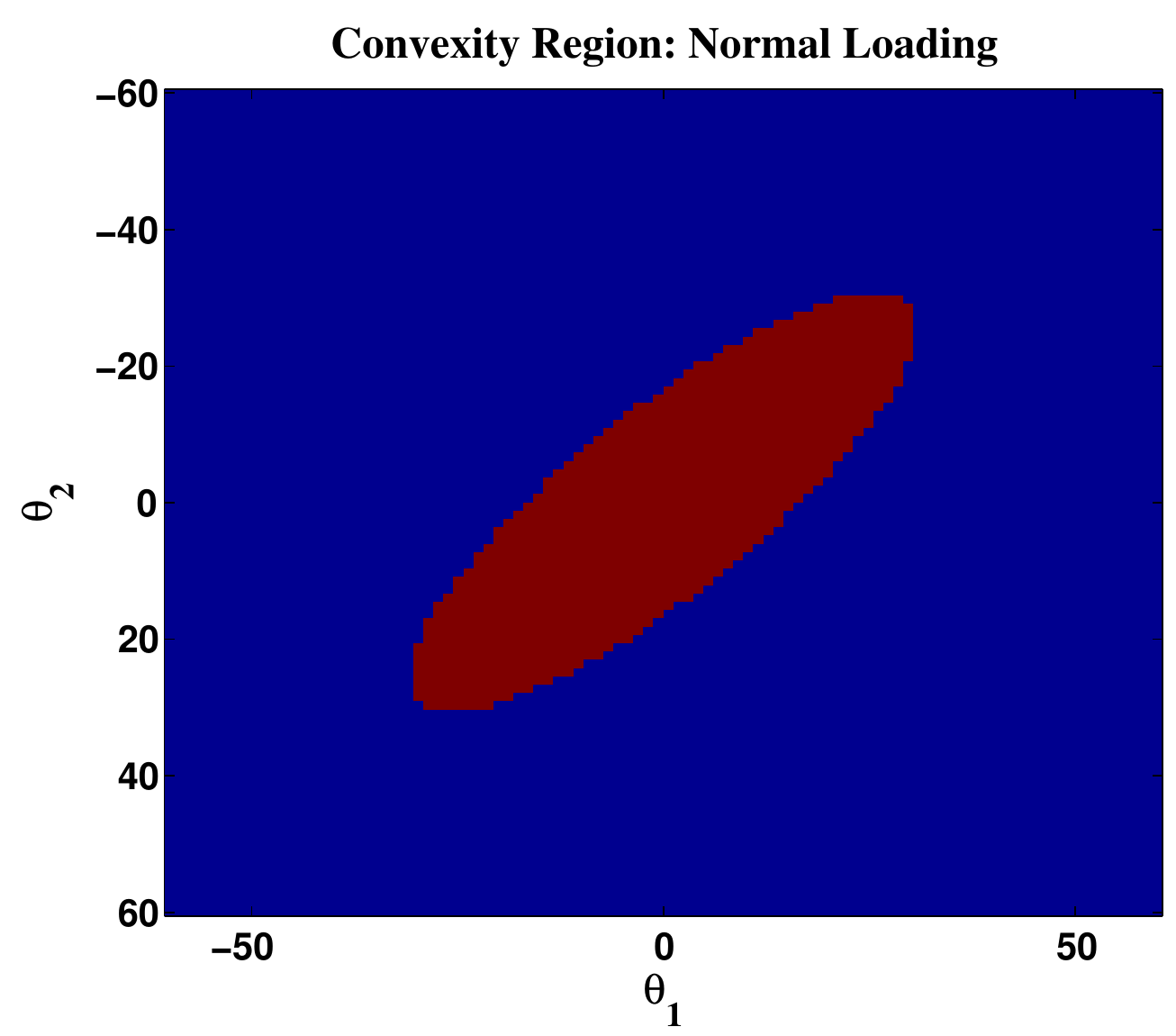}
                \caption{Regular Load}~\label{fig:ConvexityTenLoop}
        \end{subfigure}
        \begin{subfigure}[b]{0.24\columnwidth}
                \includegraphics[width=\textwidth]{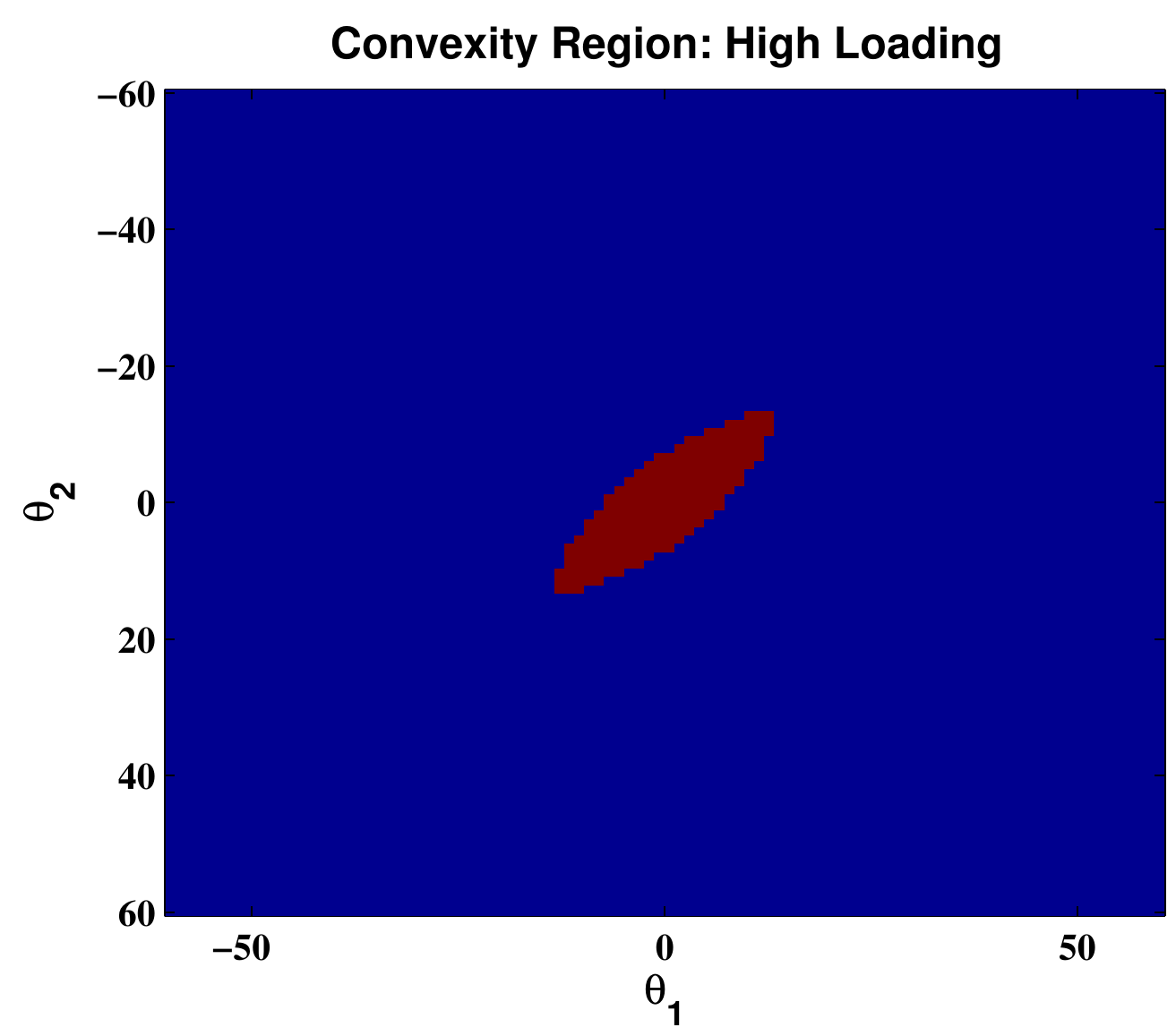}
                \caption{High Load}~\label{fig:ConvexityHunLoop}
        \end{subfigure}
         \begin{subfigure}[b]{0.24\columnwidth}
                \includegraphics[width=\textwidth]{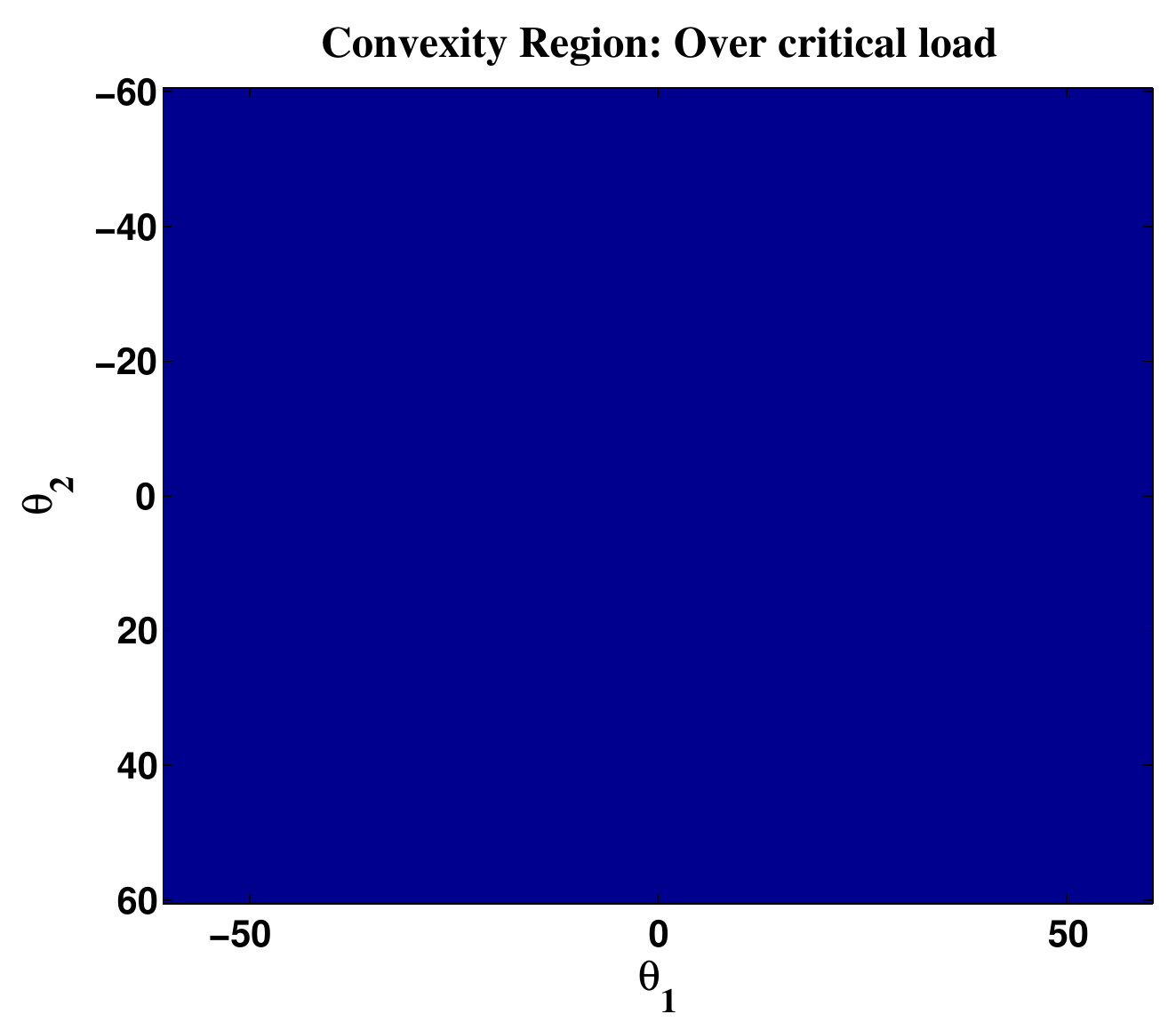}
                \caption{Loading over critical}~\label{fig:ConvexityOverLoop}
        \end{subfigure}

        \caption{Actual vs Estimated Domain of Convexity. Dark blue marks the region of non-convexity, brown represents the domain of predicted convexity (which coincides with actual convexity for all cases). One can see that as loading increases, the domain of convexity shrinks and then vanishes, indicating that there is no stable PF solution within the specified domain.}\label{fig:Convexity}
\end{figure}

Fig.~\ref{fig:Convexity} shows the actual and predicted regions of the reduced energy function convexity. For the 3-bus case considered here, our theoretical results predict the region of convexity exactly. In Fig.~\ref{fig:ConvexityNoLoop}, we plot the convexity domain for the case $B_{23}=0$. All other Figures show the case of $B_{23}>0$ with a scaled-up injection. The scaling is 3 and 5.5 and 6 in Fig.~\ref{fig:ConvexityTenLoop}, Fig.~\ref{fig:ConvexityHunLoop} and Fig.~\ref{fig:ConvexityOverLoop} respectively. As the loading increases, the region of convexity shrinks and eventually we reach the critical load where the region of convexity becomes empty. At this point, there is no longer a stable power flow solution within the given domain, since the energy function must be convex around a stable equilibrium point. These results confirm that, at least for the 3-bus case, our theorems provide an exact description for the region of convexity.

\subsection{Power Flow In Larger Networks: Existence of Solutions}

In Section \ref{sec:SolExist}, we argued that one justification for seeking solutions only within the energy function's domain of convexity is that there is an evidence to suggest that if there are no solutions within the Convexity Domain, then there are no solutions anywhere. We have a proof of this for the tree case \cite{DjMallada} and we conjecture that the result holds in general. In order to examine this hypothesis, we have studied PF solutions for two IEEE test networks: The IEEE 14 bus network and the IEEE 118 bus network. We first modify these networks to become lossless. (we simply set conductances to $0$ for all the transmission lines.) We start with the nominal base load profiles for these networks and gradually scale up the active power injections at all the buses (except the slack bus) by a factor $\kappa$ and scale up the reactive power injections by $\delta\kappa$ at all the \PQ~buses.

In order to verify insolvability of the PF equations, we use the technique developed in \cite{molzahn2012sufficient} along with the code implementing this technique in the MATPOWER package \cite{zimmerman2011matpower}. For a fixed value of $\delta$, we increase $\kappa$ until the technique from \cite{molzahn2012sufficient} detects insolvability of the power flow equation. Our convex solver does not find a PF solution when the optimal solution to \eqref{eq:optE} lies on the boundary of $\C$ or when the norm of the gradient of energy function with respect to $\br{\rho,\theta}$ is non-zero at the optimum. As $\kappa$ increases, we plot the norm of the gradient of the energy function at the optimum (a non-zero value indicates that there are no solutions in the convexity domain) and the value of the insolvability test (marking the result by $1$ in the case of insolvability  and $0$ otherwise).
\begin{figure}[htb]
        \centering
        \begin{subfigure}[b]{0.6\columnwidth}
                \includegraphics[width=\textwidth]{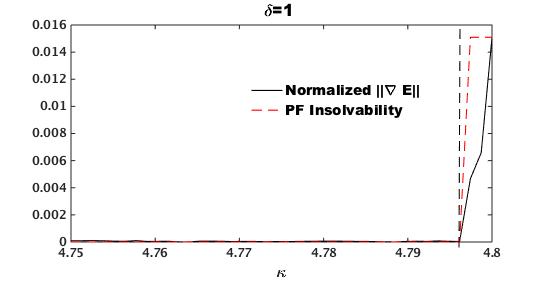}
                \caption{$\delta=1$}~\label{fig:PScaleA}
        \end{subfigure}
        \begin{subfigure}[b]{0.6\columnwidth}
                \includegraphics[width=\textwidth]{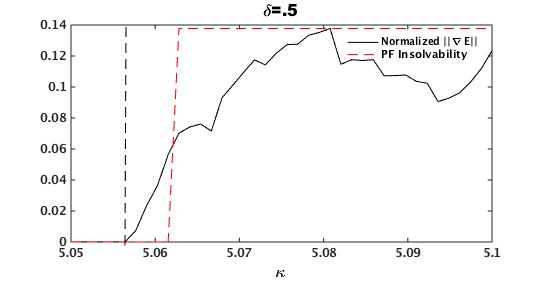}
                \caption{$\delta=.5$}~\label{fig:PScaleB}
        \end{subfigure}
        \begin{subfigure}[b]{0.6\columnwidth}
                \includegraphics[width=\textwidth]{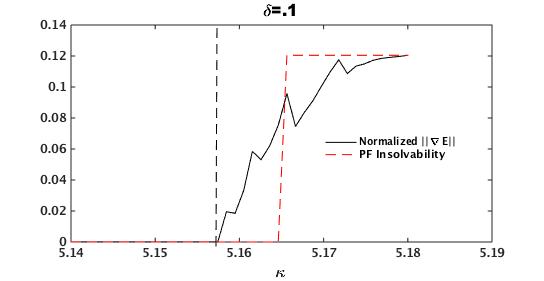}
                \caption{$\delta=.1$}~\label{fig:PScaleC}
        \end{subfigure}
        \caption{IEEE 14 Bus Case: Existence of Solutions as Load Scales. The solid black line is the norm of gradient of the energy function at the optimal solution. The value of $\kappa$ at which this becomes non-zero (marked by a dashed black line) is the point at which the convex solver fails. The red dashed line is the insolvability test. The value of $\kappa$ at which this becomes non-zero is the point at which insolvability is detected.}\label{fig:PScale14}
\end{figure}

\begin{figure}[htb]
        \centering
        \begin{subfigure}[b]{0.6\columnwidth}
                \includegraphics[width=\textwidth]{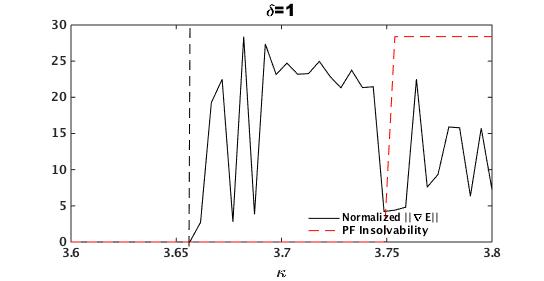}
                \caption{$\delta=1$}~\label{fig:PScale118A}
        \end{subfigure}
        \begin{subfigure}[b]{0.6\columnwidth}
                \includegraphics[width=\textwidth]{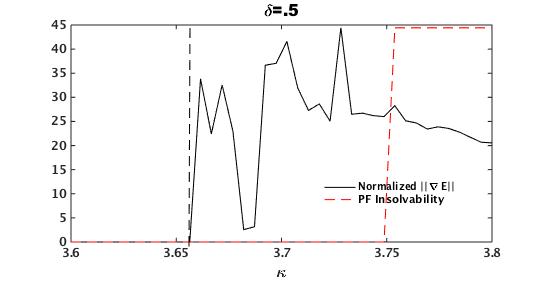}
                \caption{$\delta=.5$}~\label{fig:PScale118B}
        \end{subfigure}
        \begin{subfigure}[b]{0.6\columnwidth}
                \includegraphics[width=\textwidth]{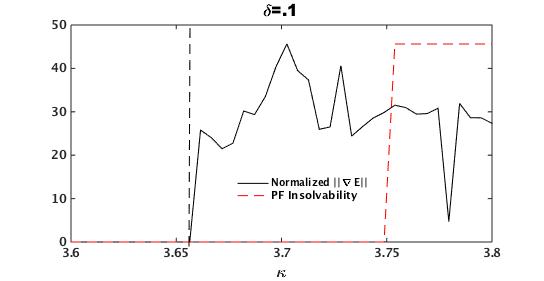}
                \caption{$\delta=.1$}~\label{fig:PScale118C}
        \end{subfigure}
        \caption{IEEE 118 Bus Case: Existence of Solutions as Load Scales.}\label{fig:PScale118}
\end{figure}

For the IEEE 14 test case, the results are shown  in  Figures \ref{fig:PScaleA},\ref{fig:PScaleB} and \ref{fig:PScaleC}. One observes that for $\delta=1$, the scaling at which our convex approach fails to find a solution is exactly the point at which the technique from \cite{molzahn2012sufficient} detects insolvability. For smaller values of $\delta$, there is a slight gap between the two, but the values are still fairly close.

For the IEEE 118 bus test case, the results are shown in Figures \ref{fig:PScaleA},\ref{fig:PScaleB} and \ref{fig:PScaleC}. We observe here (again) that the point at which \cite{molzahn2012sufficient} detects insolvability is close to the point at which our convex approach fails to find a power flow solution. The gap here is larger than the 14 bus case, but it is still rather small (difference of $.1$ in $\kappa$, which amounts to a $2.5\%$ difference in actual injections).

The results overall show that our convex approach finds solutions to the power flow equations in almost all the cases when a solution exists. There is a potential gap close to the boundary of insolvability of the PF equations. We further observed in all the cases when our algorithm reported no solution, the Newton-Raphson based approaches implemented in MATPOWER also failed to find a solution. Based on this observation, we conjecture that the gap is actually due to the approximate nature of the \cite{molzahn2012sufficient} test (based on a sufficient but not necessary guarantees for the infeasibility).

\subsection{Operational Constraints and Convexity}

The convexity condition \eqref{eq:ConvC1}\eqref{eq:ConvC2} is a complicated nonlinear matrix inequality. In this Section, we describe the simple sufficient conditions for this inequality to be satisfied by computing bounds on the phase differences $|\theta_i-\theta_j$ and voltage magnitude ratios $\expb{|\rho_i-\rho_j|}=\max\br{\frac{V_i}{V_j},\frac{V_j}{V_i}}$. We fix a bound on the log-voltage differences $\expb{|\rho_i-\rho_j|}\leq b_\rho$ and compute a bound on the phase differences $|\theta_i-\theta_j|\leq b_\theta$ so that 
\[\{\br{\rho,\theta}:\expb{|\rho_i-\rho_j|}\leq b_\rho,|\theta_i-\theta_j|\leq b_\theta\quad \forall \br{i,j}\in\E\}\subseteq \C\]
It is easy to see that the phase differences may be different for different lines, so we compute the tightest bound among all the lines.

\begin{figure}[htb]
        \centering
        \begin{subfigure}[b]{0.49\columnwidth}
                \includegraphics[width=\textwidth]{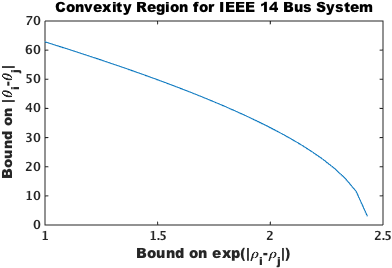}
                \caption{IEEE 14 Bus Systemß}~\label{fig:ConvexityBound14}
        \end{subfigure}
        \begin{subfigure}[b]{0.49\columnwidth}
                \includegraphics[width=\textwidth]{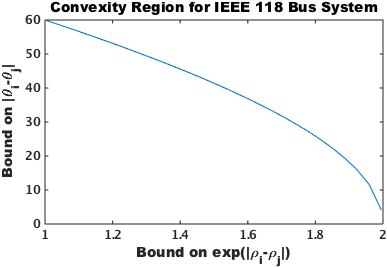}
                \caption{IEEE 118 Bus System}~\label{fig:ConvexityBound118}
        \end{subfigure}
        \caption{Convexity Dßomain and Operational Constraints}\label{fig:PScale118}
\end{figure}

For the 14 bus case, the results are plotted in Fig.~\ref{fig:ConvexityBound14}. The results show that if we allow voltage magnitude ratios across lines to be $\expb{|\rho_i-\rho_j|}\leq 1.5$, the phase differences can still be as large as $50\deg$. This is very reasonable for practical power systems, where these bounds are seldom exceeded.

For the 118 bus case, the results are plotted in Fig.~\ref{fig:ConvexityBound118}. Again, if we allow for voltage differences of up to $\expb{|\rho_i-\rho_j|}\leq 1.5$, phase differences can still be as large as $45\deg$.

To conclude, our experimental results show that with fairly lax constraints on the voltage magnitude ratios and phase differences across neighboring buses, we are guaranteed to fall within the domain of convexity. In other words, most practical solutions to the PF equations lie within the domain of convexity.

\section{Conclusions and Path Forward}\label{sec:Conclusions}

This manuscript presents novel descriptions for the domain of convexity of the energy function of a lossless power transmission network. Since the energy function provides a variational description of the power flow equations, this allows us to develop an algorithm, based on a convex optimization, solving the PF equations. This enables various applications, including detecting if the PF equations have a solution in a certain domain. There is also a growing family of important power system applications which should benefit from these convexity and PF solvability results. These applications include (a) Optimal Power Flow, (b) State and Topology Estimation, (c) Distance to Insolvability/Failure, (d) Optimal Load Shedding etc.

Another direction for future work is to extend these results to lossy networks. A special case, assuming that all lines within the network are of the same grade (the value of $X_{ij}/R_{ij}=\kappa \quad \forall\br{i,j}\in\E$), is discussed in Appendix \ref{subsec:Lossy}. In the future work we plan to exploit continuity arguments to establish existence of solutions for networks with $\frac{R}{X}$ values roughly constant. It should also be possible to quantify the degree of deviations that can be tolerated.


\section{Acknowledgments}

This work has emerged from discussions in July of 2014 at Los Alamos with Scott Backhaus and Ian Hiskens, whom the authors are thankful for guidance, comments, encouragement and criticism.  We thank Dan Molzahn and Enrique Mallada for helpful comments on various aspects of this work. We thank Florian Dorfler for pointers to useful references and comments on this manuscript. 
The work at LANL was carried out under the auspices of the National Nuclear Security Administration of the U.S. Department of Energy at Los Alamos National Laboratory under Contract No. DE-AC52-06NA25396. This material is based upon work partially supported by the National Science Foundation award \# 1128501, EECS Collaborative Research ``Power Grid Spectroscopy" under NMC. 

\bibliography{Ref}
\bibliographystyle{alpha}
\section{Appendix}
\label{sec:App}

\subsection{Proof of Theorem \ref{thm:EConvex}}
\label{subsec:Proof}

\begin{proof}
We introduce an independent edge variable $\theta_k$ for each edge $k \in \E$. Let $\theta_{\E}=\{\theta_k:k\in\E\}$. Let $k_1,k_2$ denote the ``from'' and ``to'' ends of the edge $k$. 
We then have $\theta_{\E}=A\theta$. If the energy function is convex jointly in $\br{\rho,\theta,\theta_\E}$, then it must be jointly convex in $\br{\rho,\theta}$. Writing the energy function in terms of $\br{\rho,\theta,\theta_\E}$, we derive
$E\br{\rho,\theta,\theta_{\E}}  =-\sum_{i\in\Lo\cup\G}P_i \theta_i -\sum_{i\in\Lo} Q_i\rho_i
 +\sum_{k} \frac{1}{2}B_{k}\br{\expb{2\rho_{k1}}+\expb{2\rho_{k2}}}
 -\sum_{k}B_k\expb{\rho_{k1}+\rho_{k2}}\cos\br{\theta_k}$.
In this expression only the first term depends on $\theta$. Further, the first two terms depend linearly on $\br{\rho,\theta}$ and are hence convex. Therefore, convexity of the energy function reduces to convexity of the second and third terms.  We study the Hessian which can be broken into 4 sub-matrices:
\[
\begin{pmatrix}
\nabla^2_{\rho}E\br{\rho,\theta_{\E}} & \nabla^2_{\rho,\theta_{\E}}E\br{\rho,\theta_{\E}}\\
\tranb{\nabla^2_{\rho,\theta_{\E}}E\br{\rho,\theta_{\E}}} & \nabla^2_{\theta_{\E}}E\br{\rho,\theta_{\E}}
\end{pmatrix}=\begin{pmatrix} M & N \\ \tran{N} & R\end{pmatrix}\]
where $M\in\R^{|\Lo|\times |\Lo|},N\in\R^{|\Lo|\times |\E|},R\in\R^{|\E|\times|\E|}$. For the matrix to be positive semi-definite, we require that $R\succeq 0,M-N \inve{R} \tran{N}\succeq 0 $. It is easy to see that $R$ is diagonal, with $R_{kk}=B_{k}\expb{\rho_{k1}+\rho_{k2}}\cos\br{\theta_k}$. Further,
\begin{align*}
N_{ik}&=|A_{ik}|B_k\expb{\rho_{k_1}+\rho_{k_2}}\sin\br{\theta_k},i\in\Lo,k\in\E\\
M_{ij}&=-\sum_k |A_{ik}||A_{jk}|B_{k}\expb{\rho_{k_1}+\rho_{k_2}}\cos\br{\theta_k},i\neq j\\
M_{ii}&=2B_i\expb{2\rho_i}-\sum_k |A_{ik}|B_{k}\expb{\rho_{k_1}+\rho_{k_2}}\cos\br{\theta_k}.
\end{align*}
Let $O=M-N\inve{R}\tran{N} \in \R^{|\Lo|\times|\Lo|}$
Then,
\begin{align*}
O_{ij}& =-\sum_k |A_{ik}||A_{jk}|B_{k}\expb{\rho_{k_1}+\rho_{k_2}}\br{\cos\br{\theta_k}+\frac{\sin^2\br{\theta_k}}{\cos\br{\theta_k}}}\\
&=-\sum_k |A_{ik}||A_{jk}|B_{k}\expb{\rho_{k_1}+\rho_{k_2}}\frac{1}{\cos\br{\theta_k}}
\end{align*}
for $i\neq j$ and
\begin{align*}
O_{ii}& =2B_i\expb{2\rho_i}\\
& \quad-\sum_k |A_{ik}|B_{k}\expb{\rho_{k_1}+\rho_{k_2}}\br{\cos\br{\theta_k}+\frac{\sin^2\br{\theta_k}}{\cos\br{\theta_k}}}\\
& = 2B_i\expb{2\rho_i}-\sum_k |A_{ik}|B_{k}\expb{\rho_{k_1}+\rho_{k_2}}\frac{1}{\cos\br{\theta_k}}
\end{align*}
Note also that
\begin{align*}
&L_{ij}=\frac{O_{ij}}{\expb{\rho_i}\expb{\rho_j}}=-\sum_k |A_{ik}||A_{jk}|B_{k}\frac{1}{\cos\br{\theta_k}}, i\neq j\\
&L_{ii}=\frac{O_{ii}}{\expb{2\rho_i}}=2B_i-\sum_k |A_{ik}|\expb{\rho_{k_1}+\rho_{k_2}-2\rho_i}B_{k}\frac{1}{\cos\br{\theta_k}}
\end{align*}
Then $L\succeq 0$ if and only if $O \succeq 0$ since $L=\mathrm{diag}\br{\expb{-\rho}}O\mathrm{diag}\br{\expb{-\rho}}$. Finally, note that $L\succeq 0$ if and only if the following matrix is positive semidefinite:
\begin{align*}
&\left[2B_i-\sum_{j\in\G,(i,j)\in\E}B_{ij}\expb{\rho_j-\rho_i}\frac{1}{\cos\br{\theta_{ij}}}\right]_{ii}^{\Lo}\\
& \,-\sum_{(i,j)\in\E,i,j\in\Lo}\frac{B_{ij}}{\cos\br{\theta_{ij}}} \brs{\begin{pmatrix}\expb{\rho_j-\rho_i} & 1 \\ 1 & \expb{\rho_i-\rho_j}\end{pmatrix}}_{ij}^{\Lo}
\end{align*}
To see that this is a convex constraint, we only need to observe that
\[\frac{1}{\cos\br{\theta_{ij}}}\brs{\begin{pmatrix}\expb{\rho_j-\rho_i} & 1 \\ 1 & \expb{\rho_i-\rho_j}\end{pmatrix}}_{ij}\] is $\succeq$-convex in $\br{\rho_j-\rho_i,\theta_{ij}}$ (Lemma \ref{lem:MatrixConvex}). The terms in the diagonal matrix are all concave functions of $\br{\rho,\theta}$ and hence the rerms are $\succeq$-concave. Finally, plugging in $\theta_{ij}=\theta_i-\theta_j$, one observes that the convexity is preserved and hence the energy function is convex over $\C$. Further, at $\br{\rho,\theta}=0$, the condition reduces to:
\[M\succeq 0, M_{ii}=2B_i-\sum_{j\sim i} B_{ij}=\sum_{j\sim i} B_{ij}, M_{ij}=-B_{ij}\]
Since $M$ is symmetric and diagonally dominant, it must be positive semi-definite. Hence, $0\in\C$. Further, since $\C$ is convex $\br{\alpha\rho,\alpha\theta}\in\C\quad \forall \alpha\in [0,1]$. Therefore, $\C\subset \D$.

For the converse, we look at the tree networks. We assume that the tree is connected, else the energy function can be decomposed into a sum of independent functions on each connected component and the following arguments apply. First, note that for a connected tree with $n$ buses, we find $n-1$ edges. Further, note that $\theta_{\Sb}=0$. Let $\tilde{\theta}$ denote the vector of phases at all buses except at the slack bus. Then, there is a one-to-one correspondence between $\tilde{\theta}$ and $\theta_\E$, i.e., there exists an $(n-1)\times (n-1)$ matrix $\tilde{A}$ (submatrix of $A$ formed by deleting the column corresponding to the slack bus) such that $\
\theta_\E=\tilde{A}\tilde{\theta}$. Further, note that $\tilde{A}$ is invertible. Hence, $\theta=\inv{\tilde{A}}\theta_\E$ and the Hessian of $E$ wrt $\br{\rho,\theta}$ is equal to
\begin{align*}
\begin{pmatrix}
\nabla^2_{\rho}E\br{\rho,\theta} & \nabla^2_{\rho,\theta_{\E}=\tilde{A}\theta}E\br{\rho,\theta_{\E}}\inve{\tilde{A}}\\
A^{-T}\tran{{\nabla^2_{\rho,\theta_{\E}}E\br{\rho,\theta_{\E}}}} &
{\tilde{A}}^{-T}\nabla^2_{\theta_\E=\tilde{A}\theta}E\br{\rho,\theta,\theta_\E}{\inve{\tilde{A}}}
\end{pmatrix}
\end{align*}
Using Schur-complements and the invertibility of $\tilde{A}$, it is easy to see that the positive semi-definiteness of this matrix is equivalent to that of the matrix in the first part of the proof. The positive semi-definiteness of the bottom right block requires that $|\theta_{\E}|\leq\frac{\pi}{2}$. If this inequality is strict, the block is positive definite and then the analysis in the first part is necessary and sufficient. The cases where this block is singular can be handled by a continuity argument. Hence, in this case, $\C\supseteq\D$, and combining the result from the first part, we arrive at $\D=\C$.
\end{proof}

\begin{lemma}\label{lem:MatrixConvex}
The matrix-valued function
\[f(x,y)=\frac{1}{\cos(y)}\begin{pmatrix}\expb{x} & 1 \\
1 & \expb{-x}\end{pmatrix}\]
is $\succeq$-convex, that is $
f\br{\lambda x_1+(1-\lambda)x_2,\lambda y_1+(1-\lambda)y_2}\preceq 
 \lambda f\br{x_1,y_1}+(1-\lambda)f\br{x_2,y_2}$
\end{lemma}
\begin{proof}
$f$ is $\succeq$-convex if and only if $\tr{f\br{x,y}R}$ is convex for all $R\succeq 0$ \citep{boyd2009convex}. Let $R=\begin{pmatrix}a & b \\ b & c\end{pmatrix}\succeq 0$ be an arbitrary $2\times 2$ positive semi-definite matrix. Then, we have
\begin{align*}
&\tr{f\br{x,y}R}=\frac{a\expb{x}+c\expb{-x}+2b}{\cos\br{y}}\\
& =  \frac{a\expb{x}+c\expb{-x}-2\sqrt{ac}}{\cos\br{y}}+2\frac{b+\sqrt{ac}}{\cos\br{y}}
\end{align*}
Since $R\succeq 0$, $|b|\leq\sqrt{ac}$ and hence the second term is convex (it is the inverse of a positive concave function, $\cos$). The first term can be rewritten as
$\frac{\powb{|\sqrt{a}\expb{x/2}-\sqrt{c}\expb{-x/2}|}{2}}{\cos\br{y}}$. It is easy to check that
$|\sqrt{a}\expb{x/2}-\sqrt{c}\expb{-x/2}|$ is twice differentiable at all values of $x$ and its second derivative is equal to $|\sqrt{a}\expb{x/2}-\sqrt{c}\expb{-x/2}|$. Hence, this function is convex. Therefore, the function $|\sqrt{a}\expb{x/2}-\sqrt{c}\expb{-x/2}|^2/\cos(y)$ is convex by the composition rules (since $x^2/t$ is convex in $x,t$, increasing in $x$ and decreasing in $t$). Thus, $f$ is $\succeq$-convex as long as $|y|\leq \frac{\pi}{2}$.
\end{proof}

\subsection{Lossy Networks}
\label{subsec:Lossy}

Here we show that all of our results extend to the special class of lossy networks, where all the transmission lines have the same/fixed value of the reactance-to-resistance ratio: $\frac{X_{ij}}{R_{ij}}=\frac{G_{ij}}{B_{ij}}=\kappa$. Further, we assume that all buses are \PQ~buses, except for the slack bus.

In this case, let us restate PF equations through the following linear combination of the PF equations written in their original form
\begin{subequations}
\begin{align}
& P_i+\kappa Q_i=\sum_{j \neq i} \br{\kappa^2+1}\expb{\rho_i+\rho_j}B_{ij}\sin\br{\theta_{ij}} ,i\in\Lo\label{eq:pqaLossy}\\
& \kappa P_i- Q_i=\sum_{j \neq i} \br{\kappa^2+1}\expb{\rho_i+\rho_j}B_{ij}\cos\br{\theta_{ij}}\quad \forall i \in \Lo\label{eq:pqbLossy}\\
& \theta_{\Sb}=\rho_{\Sb}=0\label{eq:sabLossy}
\end{align}\label{eq:pfLossy}
\end{subequations}
Then the PF equations \eqref{eq:pfLossy} can be rewritten in a variational form:
\begin{subequations}
\begin{align}
E\br{\rho,\theta} &=\sum_{i\in\Lo} -\br{\br{P_i+\kappa Q_i}\theta_i+\br{\kappa P_i-Q_i}\rho_i}\nonumber\\
&\quad +\sum_{\br{i,j}\in\E} B_{ij}\br{\kappa^2+1}\br{\frac{\expb{2\rho_i}+\expb{2\rho_j}}{2}}\nonumber \\
&\quad -\sum_{\br{i,j}\in\E} B_{ij}\br{\kappa^2+1}\expb{\rho_i+\rho_j}\cos\br{\theta_i-\theta_j}\label{eq:ELossy} \\
&\eqref{eq:pqaLossy}\equiv \frac{\partial E}{\partial \theta_i}=0 ,
\eqref{eq:pqbLossy}\equiv \frac{\partial E}{\partial \rho_i}=0
\end{align}
\end{subequations}

We arrive at the following statement.
\begin{theorem}
The energy function for the lossy system with constant $\frac{R}{X}$ ratio \eqref{eq:ELossy} is convex over the domain:
\begin{align}
&|\theta_i-\theta_j|\leq \frac{\pi}{2}\quad\forall (i,j)\in\E\label{eq:ConvLossyC1}\\
&\sum_{i\in\Lo} \left[2B_i-\sum_{j\in \{S\}}B_{ij}\frac{\expb{\rho_{ji}}}{\cos\br{\theta_{ij}}}\right]_{ii}^{\Lo}\nonumber \\
& -\sum_{(i,j)\in\E,i,j\in\Lo}\frac{B_{ij}}{\cos\br{\theta_{ij}}} \brs{\begin{pmatrix}\expb{\rho_{ji}} & 1 \\ 1 & \expb{\rho_{ij}}\end{pmatrix}}_{ij}^{\Lo}\succeq 0\label{eq:ConvLossyC2}
\end{align}
\end{theorem}
\begin{proof}
Almost identical to the proof of theorem \ref{thm:EConvex}.
\end{proof}

\end{document}